\documentclass[aos,seceqn,nameyear,dvips]{arximspdf}

\doi{10.1214/10-AOS796}
\volume{38}
\issue{4}
\pubyear{2010}
\firstpage{2218}
\lastpage{2241}

\makeatletter

\newtheorem{lemma}{Lemma}[section]
\newtheorem{theorem}{Theorem}[section]
\newproclaim{remark}{Remark}[section]

\newproclaim{example}{Example}

\makeatother

\begin{document}
\begin{frontmatter}

\title{Sequential monitoring of response-adaptive randomized
clinical trials}
\runtitle{Sequential monitoring of response-adaptive trials}

\begin{aug}
\author[A]{\fnms{Hongjian} \snm{Zhu}\thanksref{}\ead[label=e1]{hz5n@virginia.edu}} and
\author[A]{\fnms{Feifang} \snm{Hu}\corref{}\thanksref{t1}\ead[label=e2]{fh6e@virginia.edu}}
\runauthor{H. Zhu and F. Hu}
\affiliation{University of Virginia}
\address[A]{Department of Statistics\\
University of Virgina\\
Kerchof Hall, Charlottesville\\
Virginia 22904-4135\\
USA\\
\printead{e1}\\
\phantom{E-mail: }\printead*{e2}} %adresu isvedimo komanda gale!
\end{aug}

\thankstext{t1}{Supported by NSF Grants DMS-03-49048 and
DMS-09-07297.}

% HISTORY:
\received{\smonth{9} \syear{2009}}
\revised{\smonth{1} \syear{2010}}

% ABSTRACT
%
\begin{abstract}
Clinical trials are complex and usually involve
multiple objectives such as controlling type I error rate,
increasing power to detect treatment difference, assigning more
patients to better treatment, and more. In literature, both
response-adaptive randomization (RAR) procedures (by changing
randomization procedure sequentially) and sequential monitoring (by
changing analysis procedure sequentially) have been proposed to
achieve these objectives to some degree. In this paper, we propose
to sequentially monitor response-adaptive randomized clinical trial
and study it's properties. We prove that the sequential test
statistics of the new procedure converge to a Brownian motion in
distribution. Further, we show that the sequential test statistics
asymptotically satisfy the \textit{canonical joint distribution}
defined in Jennison and Turnbull (\citeyear{JT00}). Therefore, type I error and
other objectives can be achieved theoretically by selecting
appropriate boundaries. These results open a door to sequentially
monitor response-adaptive randomized clinical trials in practice. We
can also observe from the simulation studies that, the proposed
procedure brings together the advantages of both techniques, in
dealing with power, total sample size and total failure numbers,
while keeps the type I error. In addition, we illustrate the
characteristics of the proposed procedure by redesigning a
well-known clinical trial of maternal-infant HIV transmission.
\end{abstract}

% KEYWORDS
%

\end{frontmatter}

%s1 ###
\section{Introduction}\label{sec1}

Clinical trials usually involve multiple competing objectives such
as maximizing the power of detecting clinical difference among
treatments, minimizing total sample size and protecting more people
from possibly inferior treatments. To achieve these objectives, two
different techniques have been proposed in literature: (i) the
analysis approach---by analyzing the observed data sequentially
[sequential monitoring, Jennison and Turnbull (\citeyear{JT00})], and (ii) the
design approach---by changing the allocation probability sequentially
[response-adaptive randomization, Hu and Rosenberger (\citeyear{HR06})]. In this
paper, we discuss how to combine the two procedures in one
clinical trial in order to utilize both of their advantages.

In experiments where data accumulates sequentially, it is natural to
conduct a sequential analysis. Sequential techniques originated from a
methodology of long history based on Brownian motion.
Wald's classic work about the sequential probability ratio test
(SPRT) [Wald (\citeyear{W47})] led to the application of sequential analysis in
numerous fields of statistics. Armitage (\citeyear{A57}, \citeyear{A75}) introduced
sequential methods to clinical studies, which required monitoring
results on a patient-by-patient basis. Pocock (\citeyear{P77}) proposed
sequential monitoring of clinical trials based on a group basis.
Since then, many authors have done important work on group
sequential studies. These work are summarized in
Jennison and Turnbull (\citeyear{JT00}) and Proschan, Lan and Wittes (\citeyear{PLW06}).

The main advantages of sequential monitoring were listed in
Jennison and Turnbull (\citeyear{JT00}). First, it
is ethical to monitor clinical trials sequentially because we could
ensure that patients
are not exposed to dangerous treatments and we could stop trials as
soon as possible if needed. Second, administratively one needs to
ensure that the protocol is not violated and the assumption, which
the clinical trial is based on, is correct and valid. Third,
sequential monitoring can decrease sample size and cost.
With all the above advantages, sequential monitoring has now become a
standard technique in
conducting clinical trials.

The idea of response-adaptive randomization (RAR) can be traced back to Thompson
(\citeyear{Thompson33}) and Robbins (\citeyear{R52}). The play-the-winner rule [Zelen
(\citeyear{Z69})]
and the randomized play-the-winner rule [Wei and Durham (\citeyear{WD78})] were
proposed to reduce number of patients in the inferior treatments.
Hu and Rosenberger (\citeyear{HR03}) proved theoretically that adaptive
randomization can be used to increase statistical efficiency in some
clinical trials. In literature, many papers showed its efficient
and ethical advantages over fixed designs [Hu and Rosenberger
(\citeyear{HR06})]. With modern technology and
high capability of collecting data, it becomes easier and easier to
implement adaptive designs in sequential experiments.
Some clinical trials have already implemented the response-adaptive
designs [Rout et al. (\citeyear{Routetal93}), Tamura et al.
(\citeyear{Tamuraetal94}), Andersen (\citeyear{A96}), etc.].

Bayesian adaptive designs have also been proposed and studied in
literature. Berry (\citeyear{B05}) provided some comprehensive
introduction of Bayesian designs in clinical trials.
Recently, Cheng and Shen (\citeyear{CS05}) proposed to
sequentially monitor a Bayesian adaptive design using
decision-theoretic approaches and allowing the maximum sample size
to be sequentially adjusted by the observed data.
Lewis, Lipsky and Berry (\citeyear{LLB07}) proposed a Bayesian
decision-theoretic group sequential design for a disease with two
possible outcomes based on a quadratic loss function.
Wathen and Thall (\citeyear{WT08}) studied Bayesian adaptive model selection for
optimizing group
sequential clinical trials. In this paper, we focus on
sequential monitoring of response-adaptive randomized clinical trials.

Traditionally, sequential monitoring deals with fixed designs
(usually with equal allocation). No systematic study is available
about sequential monitoring a sequential experiment using
response-adaptive randomization, except a simulation study by Coad
and Rosenberger (\citeyear{CR99}). They found that the expected number of
treatment failures can be further reduced by combining the
triangular test with the randomized play-the-winner rule. In this
paper, we will study both theoretical properties and finite sample
properties of combining sequential monitoring with response-adaptive
randomization.

Sequential monitoring procedures use responses to stop or continue a
clinical trial. Response-adaptive randomization procedures
sequentially estimate the parameters and update the allocation
probability of the next patient. To monitor a response-adaptive
randomized clinical trial sequentially, one needs to study the two
sequential procedures simultaneously. This is conceptually difficult
because: (1) the number of patients assigned to each treatment is a
random variable at each time point; (2) both the treatment
assignments (probabilities) and the estimators of parameters (test
statistics) depend on the responses at each time point. These
problems arise from the sequential updating of estimators of the
parameters and the allocation probability function, which
leads to difficulties in finding the joint distribution of
sequential test statistics. We overcome above difficulties by (i)
approximating these different processes by martingale processes at
each time point simultaneously; (ii) then using continuous Gaussian
approximation to study these martingale processes simultaneously.

In this paper, we discuss sequential monitoring of doubly adaptive
biased coin design proposed by Hu and Zhang (\citeyear{HZ04}) for comparing two
treatments. Under widely satisfied conditions, we show that the
sequential test statistics converge to (i) a standard Brownian
motion in distribution under null hypothesis; and (ii) a~drifted
Brownian motion in distribution under alternative hypothesis. For a
standard Brownian motion, the critical value for fixed type I error
rate has been well studied in literature. Therefore,
the problem of controlling type I error is
theoretically solved. Further, we show that the sequential test
statistics satisfy the \textit{canonical joint distribution} defined in
Jennison and Turnbull (\citeyear{JT00}) asymptotically. Hence, one can
apply the group sequential methods in the book to response-adaptive
randomized clinical trials.

Simulation results support our theoretical founds in terms of type I
error and display that sequential monitoring of
response-adaptive randomization procedure could increase power and
decrease total failure number. Also compared to complete
randomization, sequential monitoring of
response-adaptive randomization procedure could stop earlier,
and thus reduce the actual sample size. In other
words, the proposed procedure achieves the goals of both RAR and
sequential monitoring. We also redesign an experiment evaluating the
effect of zidovudine treatment in reducing the risk of
maternal-infant HIV transmission performed by Connor et al. (\citeyear{Conoretal94}).
The proposed procedure can be used to decrease the number of HIV
infected people and increase the power comparing to the
complete randomization.

In Section \ref{sec2}, we introduce the notation, describe the framework and
state the main theorem. In Sections \ref{section3} and \ref{sec4}, we use both
generated data and real data to compare the proposed procedure with other
randomization procedures. Conclusions are in Section \ref{sec5} and technical
proofs are given
in the \hyperref[app]{Appendix}.

%s2 ###
\section{Sequential monitoring of response-adaptive randomization procedures}\label{sec2}

%s2.1 ###
\subsection{Notation and framework}\label{sec21}

We first describe the framework
for the randomized adaptive designs. In this article, we consider
clinical trials with two treatments 1 and 2. Let
$\mathbf{T}_i=(T_{i,1},T_{i,2})=(1,0),i=1,\ldots,n$, if the $i$th patient
is assigned to treatment 1, and $(0,1)$ otherwise, where $n$ is the
sample size.
$\mathbf{N}(n)= (N_1 (n ),N_2 (n )
)$, where
$N_j(n)=\sum_{i=1}^n T_{ij}, j=1,2$, is the number of patients in
treatment $j$. Let $\mathbf{X}=(\mathbf{X}_1,\ldots,\mathbf
{X}_n)'$, where
$\mathbf{X}_i=(\mathbf{X}_{i1},\mathbf{X}_{i2}),i=1,\ldots,n$, is a
random matrix
of responses variable and $\mathbf{X}_{ij},j=1,2$, are $d$-dimensional
random vectors. Here, only one element of $\mathbf{X}_i$, say
$\mathbf{X}_{ij}$, can be observed if the $i$th patient is assigned to
treatment $j$. We assume that $\mathbf{X}_1,\ldots,\mathbf{X}_n$ are
independent and identical distributed with unknown parameter
$(\bolds{\theta}_1, \bolds{\theta}_2)$, where $\bolds{\theta}_j$
is the
corresponding $d_j$-dimensional parameter vector
$(\theta_{j1},\ldots,\theta_{jd_j})$ of treatment $j$ ($j=1,2$). To
simplify the notation, we assume that the parameter vectors of both
treatments have the same dimension ($d_1=d_2=d$). Without loss of
generality, we also assume that $\bolds{\theta}_j=E(\mathbf{X}_{ij})$.
Otherwise, we can transform $\mathbf{X}$ and treat the transformation as
responses to make the former equation hold if such transformation
exists. Such transformation usually exists asymptotically. See
Gwise, Hu and Hu (\citeyear{GHH08}) and Hu and Zhang (\citeyear{HZ04}) for further
discussion.

Let $[nt]$ denote the largest integer that is smaller than or equal
to $nt$ for $t \in[0,1]$. Then
$\mathbf{N}([nt])= (N_1 ([nt] ),N_2 ([nt]
) )$
and $N_j([nt])=\sum_{i=1}^{[nt]} T_{ij}, j=1,2$. Note that $t=N/n$
when $N$ is the number of patients who have already been enrolled.
We introduce the so-called information time $t$ in order to
formulate this problem into the Skorohod topology [Ethier and Kurts
(\citeyear{EK86})]. After $N=[nt]$ patients have been assigned and the responses
observed, we use the modified sample means
$\hat{\bolds{\theta}}_{[nt]}=(\hat{\bolds{\theta}}_{[nt],1},\hat
{\bolds{\theta}}_{[nt],2})$
to estimate the parameter $\bolds{\theta}=(\bolds{\theta}_1,\bolds
{\theta}_2)$,
that is,
%
%e2.1 ###
\begin{equation} \label{estimators}\qquad
\hat{\bolds{\theta}}_{[nt],1}=\frac{\sum_{i=1}^{[nt]}
T_{i,1}X_{i1}+\bolds{\theta}_{0,1}}{N_1([nt])+1}   \quad\mbox{and}\quad
\hat{\bolds{\theta}}_{[nt],2}=\frac{\sum_{i=1}^{[nt]}
T_{i,2}X_{i2}+\bolds{\theta}_{0,2}}{N_2([nt])+1}.
\end{equation}
Here, we add 1 in the denominator to prevent discontinuity, and add
$\bolds{\theta}_{0,j}$, say 0.5, to estimate $\bolds{\theta}_j$
when no
patient has been assigned to the treatment $j$, $j=1,2$.

Let $\bolds{\rho}=(\rho_1,\rho_2)$ be the target allocation proportion.
Usually $\bolds{\rho}$ is obtained based on some optimal criteria and
depends on unknown parameter $\bolds{\theta}$. The selection of
$\bolds{\rho}=\bolds{\rho}(\bolds{\theta})$ has been studied by
Hayre (\citeyear{H79}),
Jennison and Turnbull (\citeyear{JT00}) and Tymofyeyev, Rosenberger and Hu
(\citeyear{TRH07}). In practice, the parameters are unknown. Therefore, we have
to first estimate them according to previous treatment assignments
and responses so that we can target the allocation proportion. We
consider a general family of doubly
adaptive biased coin design (DBCD) [Eisele and Woodroofe (\citeyear{EW95})] here.

\textit{Doubly adaptive biased coin design}: (i) assign the first
$2n_0$ patients to treatment 1 and 2 by some restricted
randomization procedures [permuted block or truncated binomial
randomization, see Rosenberger and Lachin (\citeyear{RL02})]; (ii) when the
$l$th ($l>2n_0$) patient arrives and all the responses on the
previous $l-1$ patients are available, we compute
$\hat{\bolds{\theta}}_{l-1}$ and
$\hat{\bolds{\rho}}_{l-1}=\bolds{\rho}(\hat{\bolds{\theta
}}_{l-1})$; (iii) then
assign the $l$th patient to treatment $1$ with probability
\[
g \bigl( N_1(l-1)/(l-1), \rho_1(\hat{\bolds{\theta}}_{l-1}) \bigr),
\]
where $g(s,r)\dvtx [0,1]\times[0,1]\rightarrow[0,1]$ is the allocation
function. Hu and Zhang (\citeyear{HZ04}) proposed ($\gamma\geq0$):
%
%e2.2 ###
\begin{eqnarray}\label{huzhang2004}
g^{(\gamma)}(0,r)&=&1,
\nonumber\\
g^{(\gamma)}(1,r)&=&0,
\\
g^{(\gamma)}(s,r)&=&\frac{r(r/s)^\gamma}{r(r/s)^\gamma+(1-r) (
(1-r)/(1-s) )^\gamma}.\nonumber
\end{eqnarray}
The design has drawn much attention since it was proposed and its
advantages and properties can be found in Hu and Rosenberger (\citeyear{HR03}),
Rosenberger and Hu (\citeyear{RH04}) and Tymofyeyev, Rosenberger and Hu (\citeyear{TRH07}).

To compare two treatments in clinical trials, one consider a general
hypothesis test:
\[
H_0\dvtx h(\bolds{\theta}_1)=h(\bolds{\theta}_2)
\quad\mbox{versus}\quad
H_1\dvtx h(\bolds{\theta}_1)\neq h(\bolds{\theta}_2),
\]
where $h$ is a $\Re^{d}
\rightarrow\Re$ function of parameters. In this paper, we assume
$h(\bolds{\theta}_j)$ is continuous and twice differentiable in a small
neighborhood of $\bolds{\theta}_j, j=1,2$. If one would like to test the
above hypothesis at time point $t \in(0,1]$, it is natural to
construct the test statistic as
%
%e2.3 ###
\begin{equation} \label{teststat}
Z_t \biggl(\frac{\mathbf{N}([nt])}{[nt]},\hat{\bolds{\theta
}}([nt]) \biggr)=
\frac{h (\hat{\bolds{\theta}}_1([nt]) )-h (\hat
{\bolds{\theta}}_2([nt]) )}
{\sqrt{\widehat{\operatorname{Var}} (h (\hat{\bolds{\theta}}_1([nt])
) )+
\widehat{\operatorname{Var}} (h (\hat{\bolds{\theta}}_2([nt]) ) )}}.
\end{equation}
Here,
$\widehat{\operatorname{Var}} (h (\hat{\bolds{\theta}}_1 ([nt]
)
) )$
and
$\widehat{\operatorname{Var}} (h (\hat{\bolds{\theta}}_2 ([nt]
)
) )$
are some consistent estimators of the variances of
$h (\hat{\bolds{\theta}}_1 ([nt] ) )$ and
$h (\hat{\bolds{\theta}}_2 ([nt] ) )$, respectively.
There is no covariance term on the denominator since the two terms
on the numerator are asymptotically independent [Hu, Rosenberger
and Zhang (\citeyear{HRZ06})]. Without loss of generality, we also assume that
for some functions $v_1$ and $v_2$
\[
[nt]\widehat{\operatorname{Var}} (h (\hat{\bolds{\theta}}_j ([nt]
) ) )=
v_j \biggl(\frac{\mathbf{N} ([nt] )}{[nt]},\hat{\bolds
{\theta}
} ([nt] ) \biggr)\bigl(1+o(1)\bigr)\qquad
\mbox{a.s. } j=1,2.
\]
It is easy to see that both
$v_j(\mathbf{y},\mathbf{z})$ and $Z_t(\mathbf{y},\mathbf{z})$ are
$\Re^{2+2d}\rightarrow\Re$ function, where $\mathbf{y}$ is a
two-dimensional vector and $\mathbf{z}$ is a $2d$-dimensional vector.
Examples of using this formulation are discussed in Section \ref{sec23}.

%s2.2 ###
\subsection{Main results}\label{sec22}

Based on the notation in Section \ref{sec21}, we
observe the random processes
$ (\mathbf{T}_1,\ldots,\mathbf{T}_{[nt]} )$,
$ (\mathbf{X}_1,\ldots,\mathbf{X}_{[nt]} )$, $\mathbf{N}([nt])$,
$\hat{\bolds{\theta}}_{[nt]}$, $\bolds{\rho}(\hat{\bolds{\theta
}}_{[nt]})$ and
$Z_t$ at time point $t$. When a response-adaptive randomization
procedure is used, these random processes have the following
characteristics different from those in fixed designs:
\begin{enumerate}[(2)]
\item[(1)] The allocation $ (\mathbf{N}([nt]) )$ at any
time $t$
is a random vector instead of a constant in fixed designs.
\item[(2)] The allocation $(\mathbf{N}([nt]))$ and $ (\mathbf{T}
_1,\ldots,\mathbf{T}_{[nt]} )$
are not independent with the responses
$ (\mathbf{X}_1,\ldots,\mathbf{X}_{[nt]} )$ and the parameter
estimator vector $\hat{\bolds{\theta}}_{[nt]}$.
\item[(3)] The elements $\hat{\bolds{\theta}}_1{[nt]}$ and $\hat
{\bolds{\theta}}_2{[nt]}$
depend on each other at any given time $t \in(0,1]$.
\end{enumerate}
These differences directly lead to difficulties in deriving the
joint distributions of sequential testing statistics.

To sequentially monitor a clinical trial, we need to figure out how
to control the type I error. The answer to this question relies on
the derivation of the asymptotical joint distribution of the
sequential statistics and right choices of the boundaries. Before
we give the main theorem, we need the following conditions for the
response $\mathbf{X}$, target allocation
$\rho(\bolds{\theta})$, allocation function $g$ and the function
$v_j(\mathbf{y}, \mathbf{z}),j=1,2$.

(A1) For some $\varepsilon>0$,
$E\|\mathbf{X}_1\|^{2+\varepsilon}<\infty$;

(A2) $g(s,r)$ is jointly continuous and twice differentiable at $(\rho
_1,\rho_1)$;

(A3) $g(r,r)=r$ for all $r\in(0,1)$ and $g(s,r)$ is strictly
decreasing in $s$
and strictly increasing in $r$ on $(0,1)\times(0,1)$;

(A4) $\bolds{\rho}(\mathbf{z})$ is a continuous function and twice
continuously
differentiable in a small neighborhood of $\bolds{\theta}$;

(A5) $v_j(\mathbf{y}, \mathbf{z})$ is jointly continuous and twice
differentiable in a small neighborhood of $(\bolds{\rho},\bolds
{\theta})$;

(A6) $Z_t(\mathbf{y},\mathbf{z})$ is a continuous function and it is twice
continuously differentiable in a small neighborhood of vector
$(\bolds{\rho},\bolds{\theta})$.
\begin{remark}
All the conditions are widely satisfied. An example of a design
which satisfies these conditions is DBCD in Hu and Zhang (\citeyear{HZ04}).
Condition (A1) is used to ensure the consistency of the procedure
and asymptotic normality of the allocation proportions. Condition
(A3) forces the actual allocation proportion to approach the
theoretically targeted one. Conditions (A4), (A5) and (A6) are satisfied
in all the examples in Chapter 5 of Hu and Rosenberger (\citeyear{HR06}).
\end{remark}
\begin{theorem}\label{theorem21}
Let $B_t=\sqrt{t}Z_t$ in the space $D_{[0,1]}$ with Skorohod topology.
Assume conditions \textup{(A1)}--\textup{(A6)} are satisfied. Then we have the following
two results:

\begin{longlist}
\item Under $H_0$, $B_t$ is
asymptotically a standard Brownian motion in distribution.

\item
Under $H_1$, $B_t-\sqrt{n}\mu t$ is asymptotically a standard Brownian
motion in distribution, where
\[
\mu=\frac{(h (\bolds{\theta}_1 )-h (\bolds{\theta
}_2 ))}
{\sqrt{v_1(\bolds{\rho},\bolds{\theta})+v_2(\bolds{\rho},\bolds
{\theta})}}.
\]
\end{longlist}
\end{theorem}

Based on Theorem \ref{theorem21}, we can obtain the asymptotical distribution of
the sequence of test statistics $\{Z_{t_1},\ldots,Z_{t_K}\}$, where $0
\leq t_1 \leq t_2 \leq\cdots\leq t_K \leq1$. Because
$Z_{t_i}=(\sqrt{t_i})^{-1} B_{t_i}$, we have asymptotically:

\begin{longlist}
\item
$\{Z_{t_1},\ldots,Z_{t_K}\}$ is multivariate normal;

\item $EZ_{t_i}=\mu\sqrt{nt_i}$; and

\item $\operatorname{Cov}(Z_{t_i}, Z_{t_j})=\sqrt{[nt_i]/[nt_j]}, 0\leq t_i
\leq t_j
\leq1$.
\end{longlist}

Therefore, the sequence of test statistics $\{Z_{t_1},\ldots,Z_{t_K}\}$
has the asymptotical \textit{canonical joint distribution} defined in
Jennison and Turnbull (\citeyear{JT00}).
\begin{remark}
Based on the canonical joint distribution of the sequence of test statistics
$\{Z_{t_1},\ldots,Z_{t_K}\}$, we can see that the doubly adaptive
biased coin design has a simple form of information time, which is
just the proportion of the sample size enrolled. This is because the
DBCD consistently allocates same proportion of patients to different
treatments from the beginning to the end asymptotically. We
conjecture that this simple form of information time is true for
most response-adaptive randomization procedures.
\end{remark}

Based on Theorem \ref{theorem21}, we can easily choose the correct critical
values for the asymptotic Brownian process, so that the inflation of
the type I error will be avoided. Moreover, we can also make use of
all the well-known properties of Brownian process to do further
analysis on the process of sequentially monitoring a
response-adaptive randomization procedure. Because
$\{Z_{t_1},\ldots,Z_{t_K}\}$ satisfies the \textit{canonical joint
distribution} asymptotically, we can apply the sequential techniques
in Chapters 2, 3, 4, 5, 6, 7 of Jennison and Turnbull (\citeyear{JT00}) to
response-adaptive randomized clinical trials. We may also apply
different types of spending functions to monitor a response-adaptive
randomized clinical trial sequentially. Here, we will use $\alpha$
spending functions proposed by Lan and DeMets (\citeyear{LD83}).

Any increasing function $\alpha(t)$ defined on $[0,1]$ with
$\alpha(0)=0$ and $\alpha(1)=\alpha$ is called a $\alpha$ spending
function. We spend $\alpha(t_i)-\alpha(t_{i-1})$ of the total type I
error rate at time point $t_i$, so that $\alpha(t_i)$ has been spent
after this point. For time $t_i,i=1,2,\ldots,$ we can sequentially
obtain the boundaries. This method does not require the
predetermined number of looks and equally spaced looks. We can
perform the interim monitor anytime during the trial. Such a
procedure is usually preferred by \textit{Data} and \textit{Safety}
\textit{Monitoring} \textit{Boards} (DSMB). Proschan, Lan and Wittes
(\citeyear{PLW06}) provided three special spending functions. The first one
approximates the O'Brien--Fleming boundaries [O'Brien and Fleming
(\citeyear{OF79})]
\[
\alpha_1(t)=2\{1-\Phi(z_{\alpha/2}/t^{1/2})\}.
\]
The second one is the linear spending function:
\[
\alpha_2(t)=\alpha t.
\]
The third one approximates the Pocock boundaries [Pocock
(\citeyear{P82})]:
\[
\alpha_3(t)=\alpha\ln\{1+(e-1)t\}.
\]

The O'Brien--Fleming-like function spends little of the type I error
at early looks. Consequently, the boundary for the last look is very
close to what it would have been without sequential monitoring.
Conversely, the Pocock-like function rejects the null hypothesis
easier with smaller boundaries for early looks and then has to use
a reasonably large critical value at the end to keep the type I
error. The linear function is between these two. Therefore, the
three functions above represent three typical types of spending
function. Finally, it is worth mentioning that these three spending
functions are corresponding to the process $Z_t$.

%s2.3 ###
\subsection{Examples}\label{sec23}

Here, we use two examples to illustrate how to
sequentially monitor the response-adaptive randomization procedures
based on Theorem \ref{theorem21}.
\begin{example}[(Continuous responses from normal populations)]\label{example1}
Suppose the responses of the two treatments are from
two normal distributions $Y_{i1}\sim N(\mu_1,\sigma_1^2)$ and
$Y_{i2}\sim N(\mu_2,\sigma_2^2),i=1,\ldots,n$. We would like to
compare $\mu_1$ and $\mu_2$. In this case,
$\bolds{\theta}_1=(\mu_1,\sigma_1^2+\mu_1^2)$,
$\bolds{\theta}_2=(\mu_2,\sigma_2^2+\mu_2^2)$,
$\mathbf{X}_{ij}=(Y_{ij},Y_{ij}^2)$ and
$h(\bolds{\theta}_j)=\theta_{j1}=\mu_j, j=1,2$. Then the hypotheses are
\[
H_0\dvtx \mu_1=\mu_2 \quad\mbox{versus}\quad H_1\dvtx \mu_1 \neq\mu_2.
\]
Let target allocation proportion be the Neyman allocation [Jennison and
Turnbull (\citeyear{JT00})] with
%
%e2.4 ###
\begin{equation} \label{Neymanallocation}
\rho_1=\frac{\sigma_1}{\sigma_1+\sigma_2} \quad\mbox{and}\quad
\rho_2=1-\rho_1=\frac{\sigma_2}{\sigma_1+\sigma_2}.
\end{equation}
We can use other target allocation proportions, for example,
the optimal allocation proportion [Zhang and Rosenberger
(\citeyear{ZR06})] and the $D_A$-optimal allocation proportion [Gwise, Hu and Hu
(\citeyear{GHH08})].
The sequential statistics $Z_t(\mathbf{y},\mathbf{z})$ is a function from
$\Re^{6}$ to $\Re$:
\[
Z_t(\mathbf{y},\mathbf{z})
=\frac{z_{11}-z_{21}}{\sqrt
{(z_{12}-z_{11}^2)/([nt]y_1)+(z_{22}-z_{21}^2)/([nt]y_2)}},
\]
where $\mathbf{y}={\mathbf{N}([nt])}/{[nt]}$ and
$\mathbf{z}=\hat{\bolds{\theta}}= (\hat{\theta
}_{11}([nt]),\hat
{\theta}_{12}([nt]),
\hat{\theta}_{21}([nt]),\hat{\theta}_{22}([nt]) )$. It is easy
to see that
$h (\hat{\bolds{\theta}}_1([nt]) )=\hat{\mu}_1
([nt] )$
and $h (\hat{\bolds{\theta}}_2([nt]) )=\hat{\mu
}_2([nt])$. Also
the natural variance estimators are
\[
\widehat{\operatorname{Var}} (h (\hat{\bolds{\theta}}_1 ([nt]
)
) )=\frac{\hat{\sigma}_1^2([nt])}{N_1([nt])}
\quad\mbox{and}\quad
\widehat{\operatorname{Var}} (h (\hat{\bolds{\theta}}_2 ([nt]
)
) )=\frac{\hat{\sigma}_2^2([nt])}{N_2([nt])},
\]
where $\hat{\sigma}_1^2([nt])$ and $\hat{\sigma}_2^2([nt])$ are the
usual unbiased estimators of $\sigma_1^2$ and $\sigma_2^2$ based on
the first $[nt]$ responses ($N_1([nt])$ from treatment 1 and
$N_2([nt])$ from treatment~2), respectively. Therefore,
\[
v_1(\bolds{\rho}, \bolds{\theta})=\frac{\sigma_1^2}{\rho_1}
\quad\mbox{and}\quad
v_2(\bolds{\rho}, \bolds{\theta})=\frac{\sigma_2^2}{\rho_2}.
\]
The test
statistic is then
%
%e2.5 ###
\begin{equation} \label{normaltest}
Z_t=\frac{\hat{\mu}_1([nt])-\hat{\mu}_2([nt])}
{\sqrt{\hat{\sigma}_1^2([nt])/N_1([nt])+\hat{\sigma}_2^2([nt])/N_2([nt])}}.
\end{equation}
Then based on Theorem \ref{theorem21}, the joint distribution of
$B_t=\sqrt{t}Z_t$ is asymptotically a standard Brownian process
under $H_0$. Under $H_1$, $B_t-\sqrt{n}\mu t$ is asymptotically a
standard Brownian motion in distribution, where
\[
\mu=\frac{\mu_1-\mu_2}
{\sqrt{\sigma_1^2/\rho_1+\sigma_2^2/(1-\rho_1)}}.
\]
\end{example}
\begin{example}[(Binary responses)]\label{example2}
Assume $Y_{i1} \sim \operatorname{Bin}(1,
p_1)$ and
$Y_{i2} \sim \operatorname{Bin}(1$,\break $p_2),i=1,\ldots,n$, and we would like to compare
$p_1$ and $p_2$. In this case,
$\bolds{\theta}_1=(p_1)$, $\bolds{\theta}_2=(p_2)$, $\mathbf
{X}_{ij}=(Y_{ij})$
and $h(\bolds{\theta}_j)=\theta_{j1}, j=1,2$. The hypotheses are
\[
H_0\dvtx p_1=p_2 \quad\mbox{versus}\quad H_1\dvtx p_1 \neq p_2.
\]
Three common target allocations are: (i) Neyman allocation,
%
%e2.6 ###
\begin{eqnarray} \label{Neymanbin}
\rho_1&=&\frac{\sqrt{p_1(1-p_1)}}{\sqrt{p_1(1-p_1)}+\sqrt
{p_2(1-p_2)}}\qquad
\mbox{and}\nonumber\\[-8pt]\\[-8pt]
\rho_2&=&\frac{\sqrt{p_2(1-p_2)}}{\sqrt{p_1(1-p_1)}+\sqrt{p_2(1-p_2)}};\nonumber
\end{eqnarray}

\mbox{}\phantom{i}(ii) optimal allocation proposed by Rosenberger et al. (\citeyear{Rosenbergeretal01}),
%
%e2.7 ###
\begin{equation} \label{Rosopt}
\rho_1=\frac{\sqrt{p_1}}{\sqrt{p_1}+\sqrt{p_2}}  \quad\mbox{and}\quad
\rho_2=\frac{\sqrt{p_2}}{\sqrt{p_1}+\sqrt{p_2}};
\end{equation}

(iii) Urn allocation [Wei and Durham (\citeyear{WD78})],
%
%e2.8 ###
\begin{equation} \label{urnallocation}
\rho_1=\frac{q_2}{q_1+q_2}  \quad\mbox{and}\quad
\rho_2=\frac{q_1}{q_1+q_2}.
\end{equation}

Neyman allocation is a commonly discussed allocation which is related
to the efficiency issue in the field of response-adaptive
randomization procedures. We studied sequential monitoring of
response-adaptive designs with Neyman allocation in order to show
that our proposed procedure is able to achieve various objects.

In this case, $Z_t(\mathbf{y},\mathbf{z})$ is a function from $\Re
^{4}$ to
$\Re$:
\[
Z_t(\mathbf{y},\mathbf{z})=\frac{z_{11}-z_{21}}
{\sqrt{z_{11}(1-z_{11})/([nt]y_1)+z_{21}(1-z_{21})/([nt]y_2)}},
\]
where $\mathbf{y}=(N_1([nt])/[nt],N_2([nt])/[nt])$,
$\mathbf{z}= (\hat{\theta}_{11}([nt]),\hat{\theta
}_{21}([nt]) )$,
$h (\hat{\bolds{\theta}}_1([n\times\break t]) )=\hat{p}_1([nt])$ and
$h (\hat{\bolds{\theta}}_2([nt]) )=\hat{p}_2([nt])$. The
corresponding variance estimators are
\[
\widehat{\operatorname{Var}} (h (\hat{\bolds{\theta}}_1 ([nt]
)
) )
=\frac{\hat{p}_1 ([nt] ) (1-\hat{p}_1
([nt] ) )}{N_1 ([nt] )}
\]
and
\[
\widehat{\operatorname{Var}} (h (\hat{\bolds{\theta}}_2 ([nt]
)
) )
=\frac{\hat{p}_2 ([nt] ) (1-\hat{p}_2([nt])
)}{N_2([nt])}.
\]
Therefore,
\[
v_1(\bolds{\rho}, \bolds{\theta})=\frac{p_1(1-p_1)}{\rho_1}
\quad\mbox{and}\quad
v_2(\bolds{\rho}, \bolds{\theta})=\frac{p_2(1-p_2)}{\rho_2}.
\]
The test statistic is
%
%e2.9 ###
\begin{eqnarray} \label{bintest}
Z_t&=&\bigl({\hat{p}_1([nt])-\hat{p}_2([nt])}\bigr)\nonumber\\[-8pt]\\[-8pt]
&&{}\times\Biggl(\sqrt{\hat{p}_1([nt])
\frac{1-\hat{p}_1 ([nt] ) }{N_1([nt])}+\hat
{p}_2([nt]) \frac{1-\hat{p}_2([nt]) }{N_2([nt])}}\Biggr)^{-1}.\nonumber
\end{eqnarray}
Then $B_t=\sqrt{t}Z_t$ converges to a standard Brownian process in
distribution under~$H_0$. Under $H_1$,
$B_t-\sqrt{n}\mu t$ is asymptotically a standard Brownian motion in
distribution, where
\[
\mu=\frac{p_1-p_2}
{\sqrt{p_1(1-p_1)/\rho_1+p_2(1-p_2)/(1-\rho_1)}}.
\]

Theorem \ref{theorem21} can be applied to different situations such as the
examples considered in Chapter 5 of Hu and Rosenberger (\citeyear{HR06}). In
Examples \ref{example1} and \ref{example2}, now assume we would like to look at the process at
three points: $t_1=0.2$, $t_2=0.5$ and $t_3=1$. Then we can use the
corresponding critical values from the three spending functions
[Proschan, Lan and Wittes (\citeyear{PLW06})] in the last subsection for $Z_t$ to
keep the overall type I error 0.05: O'Brien--Fleming-like boundaries
(4.877, 2.963, 1.969), linear boundaries (2.576, 2.377, 2.141) and
Pocock-like boundaries (2.438, 2.333, 2.225).
\end{example}

%s3 ###
\section{Simulation study}\label{section3}

In Section \ref{sec2}, we obtained the asymptotical distribution of the test
statistic $Z_t$. In this section, we will use the two examples in
Section \ref{sec2} to study the finite sample properties of the proposed procedure.

In Examples \ref{example1} and \ref{example2}, we use the doubly adaptive biased coin design
with Hu and Zhang's allocation function in (\ref{huzhang2004}) and
$\gamma=2$ is used. In Tables \ref{table1}--\ref{table5}, we use the same total sample size
500. The first $50$ patients ($n_0=25$) are randomly assigned to
treatments 1 and 2 by using permuted block randomization. Then, for the
$l$th ($l>50$) patient, the unknown parameters are estimated by
using (\ref{estimators}) based on the first $l-1$ responses with
$\theta_{0,1}=\theta_{0,2}=0.5$. For normal responses in Example~\ref{example1},
we estimate $\sigma_1^2$ and $\sigma_2^2$ by using the standard
unbiased estimators based on the first $l-1$ responses.

%t1 ###
\begin{table}
\caption{Example \protect\ref{example1} with Neyman allocation, $\mu_1=\mu_2=1$, $\sigma
_1=1$, $\sigma_2=2$}\label{table1}
\begin{tabular*}{\tablewidth}{@{\extracolsep{\fill}}lccc@{}}
\hline
\textbf{Critical values} & \textbf{Randomization} &
\textbf{Type I error} & $\bolds{\hat{\rho}_1}$ \textbf{(s.e.)}\\
\hline
B--F-like & DBCD &0.055 &0.333 (0.020)\\
B--F-like & CR&0.052 & 0.500 (0.022) \\
Linear & DBCD &0.048 &0.333 (0.020)\\
Linear & CR &0.053 &0.500 (0.023) \\
Pocock-like & DBCD &0.051& 0.332 (0.020) \\
Pocock-like & CR &0.052&0.500 (0.023) \\
\hline
\end{tabular*}
\end{table}

%t2 ###
\begin{table}[b]
\caption{Example \protect\ref{example2} with optimal allocation, $p_1=p_2=0.5$}\label{table2}
\begin{tabular*}{\tablewidth}{@{\extracolsep{\fill}}lccc@{}}
\hline
\textbf{Critical values} & \textbf{Randomization} & \textbf{Type I error}
& $\bolds{\hat{\rho}_1}$ \textbf{(s.e.)}
\\
\hline
B--F-like & DBCD &0.051 &0.500 (0.016)\\
B--F-like & CR&0.046 & 0.500 (0.023) \\
Linear & DBCD &0.055 &0.500 (0.019) \\
Linear & CR &0.061 &0.500 (0.023) \\
Pocock-like & DBCD &0.056& 0.500 (0.019) \\
Pocock-like & CR &0.050&0.500 (0.022) \\
\hline
\end{tabular*}
\end{table}

%t3 ###
\begin{table}
\caption{Example \protect\ref{example1} with Neyman allocation, $\mu_1=1$, $\mu_2=1.4$,
$\sigma_1=1$, $\sigma_2=2$}\label{table3}
\begin{tabular*}{\tablewidth}{@{\extracolsep{\fill}}lcccccc@{}}
\hline
\textbf{Critical values} & \textbf{Randomization} & \textbf{Power}
& $\bolds{\hat{\rho}_1}$ \textbf{(s.e.)}
& $\bolds{N_1}$ & $\bolds{N_2}$ & $\bolds{N_3}$
\\
\hline
B--F-like & DBCD &0.847 &0.333 (0.021)&\phantom{00}2 & 1013 & 3222\\
B--F-like & CR&0.807 & 0.500 (0.024) & \phantom{00}1 &\phantom{0}842 &3193\\
Linear & DBCD &0.812 &0.332 (0.027)&594 & 1429 & 2035\\
Linear & CR &0.765 &0.500 (0.028) &477 & 1380 & 1970\\
Pocock-like & DBCD &0.792& 0.332 (0.028) &741 & 1443 &1774\\
Pocock-like & CR &0.738&0.500 (0.028) &544 & 1309 & 1835\\
\hline
\end{tabular*}
\end{table}

%t4 ###
\begin{table}[b]
\tabcolsep=0pt
\caption{Example \protect\ref{example2} with urn allocation, $p_1=0.5$, $p_2=0.625$}\label{table4}
\begin{tabular*}{\tablewidth}{@{\extracolsep{\fill}}lccccccc@{}}\hline
\textbf{Critical values} & \textbf{Randomization} & \textbf{Power} & $\bolds{\hat{\rho}_1}$ \textbf{(s.e.)}
& $\bolds{N_1}$& $\bolds{N_2}$ & $\bolds{N_3}$ & \textbf{Total failures
(s.e.)}
\\
\hline
B--F-like & DBCD &0.811&0.426 (0.033)&\phantom{00}4 & 839 & 3214 & 211 (13)\\
B--F-like & CR&0.811 & 0.500 (0.024)& \phantom{00}1 &839 &3215& 217 (13)\\
Linear & DBCD &0.762 &0.421 (0.041)&503 & 1396 & 1912& 206 (14)\\
Linear & CR &0.767 &0.500 (0.029) &521 & 1300 & 2016 &212 (14)\\
Pocock-like & DBCD &0.749& 0.421 (0.042) &609 & 1325 &1809& 205 (14)\\
Pocock-like & CR &0.738&0.501 (0.029) &603 & 1312 & 1773 & 211 (15)\\
\hline
\end{tabular*}
\end{table}

For simplicity, we look at the test at three time points [$n_1=100$
$(t_1=0.2)$, $n_2=250$ $(t_2=0.5)$ and $n=500$ $(t_3=1)$]. Then the
three sets of spending function boundaries in Section \ref{sec23} are used
to ensure $\alpha= 0.05$. For each spending function, the first row
in the table is for DBCD and the second row is for complete
randomization (denoted as CR in the tables). All the simulations are
based on 5000 replications.

In Table \ref{table1}, we simulate Example \ref{example1} with two normal responses $N(1,1)$
and $N(1,2)$ by using the Neyman allocation
(\ref{Neymanallocation}). We find that the type I error of
sequentially monitoring
the response-adaptive randomization procedure and complete
randomization are both well kept at the 0.05 level. We also report
the mean and standard deviation of actual allocation proportion
($\hat{\rho}_1$) for treatment 1 [$N(1,1)$]. We find that the mean
agrees with Neyman allocation and the standard deviation is
reasonably small for DBCD. This indicates that the DBCD is able to
target the theoretical targeted allocation proportion very well. In
Table \ref{table2}, we simulate the Example \ref{example2} with two binary responses
$p_1=p_2=0.5$ and the target allocation is the optimal allocation
(\ref{Rosopt}). We obtain the same conclusion as Table \ref{table1}. We have
also done simulations for some other cases, and similar results are
obtained. These numerical results indicate that sequential
monitoring of the response-adaptive randomization will not inflate
the type I error with the appropriate boundaries based on Theorem
\ref{theorem21}.

%t5 ###
\begin{table}
\tabcolsep=0pt
\caption{Example \protect\ref{example2} with optimal allocation, $p_1=0.5$, $p_2=0.625$}\label{table5}
\begin{tabular*}{\tablewidth}{@{\extracolsep{\fill}}lccccccc@{}}
\hline
\textbf{Critical values} & \textbf{Randomization} & \textbf{Power}
& $\bolds{\hat{\rho}_1}$ \textbf{(s.e.)} &
$\bolds{N_1}$ & $\bolds{N_2}$ & $\bolds{N_3}$ & \textbf{Total failures
(s.e.)}
\\
\hline
B--F-like & DBCD &0.810&0.471 (0.017) &\phantom{0}4 & \phantom{0}863 & 3185 & 214 (12)\\
B--F-like & CR&0.805 & 0.501 (0.024) & \phantom{0}4 &\phantom{0}795 &3229& 218 (13)\\
Linear & DBCD &0.768 &0.468 (0.022) &520 & 1354 & 1964& 210 (14)\\
Linear & CR &0.762 &0.500 (0.029) &474 & 1367 & 1971 & 214 (14)\\
Pocock-like & DBCD &0.754& 0.469 (0.023) &673 & 1309 &1787& 210 (14)\\
Pocock-like & CR &0.749&0.500 (0.03)\phantom{0}& 602 & 1351 & 1793 & 213 (15)\\
1.96 & DBCD & 0.805 & 0.472 (0.015)&NA&NA&NA& 217 (11)\\
1.96 & CR & 0.802 & 0.500 (0.022)&NA&NA&NA& 221 (11)\\
\hline
\end{tabular*}
\end{table}

Next, we show other advantages of the sequential monitoring of the
response-adaptive
randomization procedure.
In Table \ref{table3}, we simulate Example \ref{example1} with two normal responses $N(1,1)$
and $N(1.4,2)$ using Neyman allocation (\ref{Neymanallocation}) as the
target allocation that
maximizes the power. The power of the sequential monitoring of the
response-adaptive randomization procedure is about 5\%--8\% higher
than sequentially monitoring the complete randomization. $N_i$ in the
table is the number of rejections at the $i$th look. Rejection at the
first two looks means stopping the trial earlier. DBCD with
sequential monitoring obviously stops the trial earlier than
complete randomization.

In Table \ref{table4}, we simulate Example \ref{example2} with two binary responses
$p_1=0.5$ and $p_2=0.625$ using the urn allocation
(\ref{urnallocation}) as the target allocation that assigns more
people to the better treatment. If we reject the null hypothesis at
the first two looks, we assign all the remaining patients to the
estimated better treatment and count the total failure number. We do
this only for the comparison in the simulation study. In a real
clinical trial, we stop the trial if the null hypothesis is rejected
at an interim look. From the mean total failure number, the
DBCD with sequential monitoring has lower failure numbers than
complete randomization for each type of spending function. $N_1$,
$N_2$, and $N_3$ show that our methods stop the trial a little
earlier and the power is almost the same.

In Table \ref{table5}, we simulate Example \ref{example2} with two binary responses
$p_1=0.5$ and $p_2=0.625$ using the optimal
allocation (\ref{Rosopt}) used to maximize the power while keeping the
total failure number. We deal with the remaining patients in the
same way as in Table \ref{table4} if we reject the null hypothesis at the first
two looks. We find that sequential monitoring of the response-adaptive
randomization procedure can achieve the aim of optimal allocation.
Its power is larger and its failure number is less than the complete
randomization procedure. In this table, we also do the simulation
without sequential monitoring. That is, we only look at the test
once at the end of the trial and the critical value is 1.96 for the
nominal significance level 0.05. We report it at the last two rows.
It is obvious that sequential monitoring can reduce the total
failures.

Based on the simulation results, we can see the
advantages of sequentially monitoring response-adaptive
randomized clinical trials: (i) controlling type I error well;
(ii) reducing the total number of failures; (iii) increasing power; and
(iv)~stopping the trail earlier (reducing total sample size).

%s4 ###
\section{Re-designing the HIV transmission trial}\label{sec4}

Maternal-infant transmission is the primary means by which infants
are infected by HIV virus. Connor et al. (\citeyear{Conoretal94}) reported a trial to
evaluate the drug AZT (Zidovudine treatment) in reducing the risk of
maternal-infant HIV transmission. In this clinical trial, 477
HIV-infected pregnant women were enrolled from April 1991 to
December 1993 and assigned to the Zidovudine treatment group and
placebo group with a 50--50 randomization scheme. This experiment was
a randomized, double-blind and placebo-controlled trial. 239 were
allocated to the treatment group and 238 to the placebo group. At
the end of the trial, $8.3\%$ of the infant from the treatment group
were infected by the HIV virus, while $25.5\%$ from the placebo
group were infected.

%t6 ###
\begin{table}
\caption{Re-designed the HIV trial with full sample size}\label{table6}
\begin{tabular*}{\tablewidth}{@{\extracolsep{\fill}}lcccc@{}}
\hline
\textbf{Target allocation} & \textbf{Critical values}
& $\bolds{\hat{\rho}_1}$ \textbf{(s.e.)}
& \textbf{Power} & \textbf{Total failures (s.e.)}
\\
\hline
CR&linear & 0.500 (0.039) & 0.999 & 60.1 (11.1)\\
CR&1.96 & 0.501 (0.023)&0.999&80.7 (8.2)\phantom{0}\\
Urn allocation& linear & 0.751 (0.062)& 0.996 & 52.3 (9.2)\phantom{0}\\
Optimal allocation &linear & 0.527 (0.021)& 0.997& 56.4 (10.8)\\
\hline
\end{tabular*}
\end{table}

In Table \ref{table6}, we redesign the study by sequential monitoring of both
complete randomization (the first two rows in the table) and
response-adaptive randomization [DBCD (\ref{huzhang2004}) with $\gamma
=2$] (the last
three rows in the table). We assume the success rate for the
treatment group is $p_1=0.917$ and that for the placebo group is
$p_2=0.745$ (as reported in the original paper). We look at the test
at the three same time points as mentioned in the last section,
$n_1=95$ $(t_1=0.2)$,
$n_2=143$ $(t_2=0.5)$ and $n=239$ $(t_3=1)$. The boundary we use is the
linear spending function (2.576, 2.377, 2.141) except the second row
in the table where we do the equal allocation without sequential
monitoring. We report the actual allocation proportion for the
treatment group,
power and the total HIV-infected number. As before, if we reject the
null hypothesis at the first two looks, we will assign all the
remaining patients to the estimated better treatment. We find that the
sequential
monitoring technique will decrease the HIV-infected number
dramatically from the first two rows. Response-adaptive
randomization technique will also reduce the HIV-infected number
compared to the complete randomization. Sequential monitoring DBCD
while targeting at the urn allocation has the least HIV-infected number,
which agrees with the aim of urn allocation.

%t7 ###
\begin{table}
\caption{Re-designed the HIV trial with sample size $n=245$}\label{table7}
\begin{tabular*}{\tablewidth}{@{\extracolsep{\fill}}lcccc@{}}
\hline
\textbf{Target allocation} & \textbf{Critical values} & $\bolds{\hat{\rho}_1}$
\textbf{(s.e.)} & \textbf{Power} & \textbf{Total failures (s.e.)}
\\
\hline
CR & B--F-like & 0.500 (0.036) &0.947 & 40.1 (7.0)\\
CR & linear & 0.501 (0.042) &0.942 & 36.6 (7.5)\\
CR & 1.96 & 0.500 (0.032) &0.958& 43.1 (5.8)\\
Urn allocation & B--F-like & 0.745 (0.068) &0.920& 30.7 (5.9)\\
Urn allocation & linear & 0.747 (0.074) &0.885&29.3 (6.1)\\
Optimal allocation & B--F-like & 0.528 (0.023)&0.952& 36.8 (6.7)\\
Optimal allocation & linear & 0.529 (0.025)&0.945& 32.8 (7.3)\\
\hline
\end{tabular*}
\end{table}

In Table \ref{table7}, we reduce the full sample size to 245
(to achieve power $0.95$ for complete randomization) and keep all the
other settings unchanged. We obtain the same conclusion about the HIV-infected
number as in Table \ref{table6}. We also find that targeting optimal allocation
with DBCD has slightly higher power than targeting equal allocation
when sequential monitoring is used. Targeting urn allocation with
DBCD has slightly less power but the HIV-infected number in this way
is the least. Overall, sequential monitoring of the
response-adaptive randomization procedure is better than that of
complete randomization, since it reduces the HIV-infected
number and remains good power.

%s5 ###
\section{Conclusion remarks}\label{sec5}

Now sequential monitoring becomes a standard technique in clinical
trials. To apply response-adaptive randomization in clinical trials,
it is important to know how to sequentially monitor adaptive
randomized trials. In this paper, we overcome this hurdle and show
the advantages of sequential monitoring response-adaptive randomized
clinical trials both theoretically and numerically. We use a Gaussian
process in the Skorohod topology to describe the relationship
between the allocation and parameter estimators. One of the main
contributions of this paper is to show that sequential statistics can
be asymptotically approximated by a Brownian process in distribution
under both null and alternative hypotheses. Further, we find that
the sequential test statistics satisfy the \textit{canonical joint
distribution} asymptotically. Consequently, the results of this
paper not only solve the problem of preserving a preset type I
error but may lead to many area of potential future research.

We have studied how to sequentially monitor a clinical trial
based on doubly adaptive biased coin design proposed by Eisele and
Woodroofe (\citeyear{EW95}) and Hu and Zhang (\citeyear{HZ04}). Another important family
of response-adaptive randomization procedure is based on urn models,
which include randomized play-the-winner rule [Wei and Durham
(\citeyear{WD78})], generalized Friedman's urn models [Athreya and Karlin (\citeyear{AK68}), Bai
and Hu (\citeyear{BH05})],
drop-the-loser rule [Ivanova (\citeyear{I03})], sequential estimation-adjusted
urn models [Zhang, Hu and Cheung (\citeyear{ZHC06})], etc. The technique used in
this paper opens a door to study the properties of
sequential monitoring of clinical trials based on
these urn models or the efficient randomized adaptive designs [Hu,
Zhang and He (\citeyear{HZH09})]. We leave this for future study.

In this paper, we have used $\alpha$-spending function to calculate
the critical boundaries. Because the sequential test statistics
satisfy the \textit{canonical joint distribution} asymptotically, we
can implement all the sequential techniques introduced in Jennison
and Turnbull (\citeyear{JT00}) based on this canonical form. Also we can use
the optimal spending functions in Anderson (\citeyear{A07}), or the beta spending
functions in DeMets (\citeyear{DeMets06}). We also leave the details for future
research.

%%%%%%%%%%%%%%%%%%%%%%%%%%%%%%%%%%%%%%%%%%%%%%%%%%%%%%%%%%%%%%%%%%%%%%%%%%%%%%%%%%%%%%%%%%%%%%
\begin{appendix}\label{app}
\section*{Appendix: Proofs}

First, we introduce some further notation. For a function
$\bolds{\eta}(\mathbf{u},\mathbf{w})\dvtx\Re^L \times\Re^M
\rightarrow
\Re^2$, we
denote the partial derivative matrices as
\[
\nabla_u(\bolds{\eta})= \biggl(\frac{\partial\eta_k}{\partial u_i};
i=1,\ldots,L, k= 1,2 \biggr)_{L\times2}
\]
and
\[
\nabla_w(\bolds{\eta})= \biggl(\frac{\partial\eta_k}{\partial w_j};
j=1,\ldots,M,
k= 1,2 \biggr)_{M\times2}.
\]
Let
$H=\nabla_r (g(r,s),1-g(r,s) )|_{(\rho_1,\rho_1)}$ and
$E=\nabla_s(g(r,s),1-g(r,s))|_{(\rho_1,\rho_1)}$ be the partial
derivative matrices of the allocation function $g$. Further, let
$V=\operatorname{diag} (\operatorname{var}(\mathbf{X}_{11})/{\rho_1},\operatorname{var}(\mathbf
{X}_{12})/{\rho
_2} )$,
$\Sigma_3= (\nabla(\bolds{\rho})|_{\bolds{\theta}}
)'V\nabla
(\bolds{\rho})|_{\bolds{\theta}}$, $\Sigma_1=\operatorname{diag}(\bolds{\rho
})-\bolds{\rho}
'\bolds{\rho}$
and $\Sigma_2=E'\Sigma_3E$. In Hu and Zhang (\citeyear{HZ04}), they studied the
asymptotic properties of $\mathbf{N}(n)$, $\hat{\bolds{\rho}}(n)$ and
$\hat{\bolds{\theta}}(n)$ at the end of the trial. Based on their
results, one can do the corresponding statistical inference after
observing all responses of the clinical trial. To monitor the
response-adaptive randomized trial sequentially, we need to know the
theoretical properties of the process $\mathbf{N}([nt])$ and
$\hat{\bolds{\theta}}([nt])$ for any given $t \in(0,1]$. To do
this, we start
with Lemma \ref{lemma1}.
\begin{lemma}\label{lemma1}
Let $W_{1t}$ and $W_{2t}$ be two independent standard\break
two-dimensional Brownian processes. $\mathbf{N}([nt])$,
$\hat{\bolds{\theta}}([nt])$, $\bolds{\rho}$ and $\bolds{\theta
}$ are defined as
in Section \ref{sec2}. Under the conditions of Theorem \ref{theorem21}, we have
%
%e5.1 ###
\begin{equation}
n^{-1/2}([nt]) \biggl(\frac{\mathbf{N}([nt])}{[nt]}-\bolds{\rho},
\hat{\bolds{\theta}}([nt])-\bolds{\theta} \biggr) \rightarrow
(G_t,W_{2t}V^{1/2})
\end{equation}
in distribution in the space $D_{[0,1]}$ with the Skorohod topology,
where the Gaussian process
%
%e5.2 ###
\begin{equation}
G_t=\int_0^t (dW_{1x} )\Sigma_1^{1/2} \biggl(\frac
{t}{x} \biggr)^H+
\int_0^t (dW_{2x} )\Sigma_2^{1/2}\biggl[\int_x^t\frac
{1}{y} \biggl(\frac{t}{y} \biggr)^H\,
dy\biggr],
\end{equation}
which is the solution of the stochastic differential equation
\[
dG_t=(dW_{1t}) \Sigma_1^{1/2} + \frac{W_{2t}\Sigma_2^{1/2}}{t} \,dt+
\frac{G_t}{t} H \,dt \qquad\mbox{with } G_0=0,
\]
and $a^H$ is the matrix power function defined as
\[
a^H=e^{H\ln a}=\sum_{j=0}^\infty\frac{(\ln a)^j}{j!}H^k.
\]
\end{lemma}
\begin{pf}
It is worth noting that the
response-adaptive design in Theorem~\ref{theorem21} satisfies all the conditions
of Hu and Zhang (\citeyear{HZ04}). So all the results in Hu and Zhang (\citeyear{HZ04})
are valid. We will prove this lemma by using the weak convergence of
the martingale [cf. Theorem 4.1 of Hall and Heyde (\citeyear{HH80})]. To do
this, we first approximate the process
$ (\frac{\mathbf{N}([nt])}{[nt]}-\bolds{\rho},
\hat{\bolds{\theta}}([nt])-\bolds{\theta} )$ by a\vspace*{1pt} martingale
and then
prove the following two facts: (1) Lindeberg condition holds for the
approximated martingale process; and (2) the limiting covariance
of $n^{-1/2}([nt]) (([nt])^{-1}\mathbf{N}([nt])-\bolds{\rho},
\hat{\bolds{\theta}}([nt])-\bolds{\theta} )$ agrees with
that of
$(G_t,W_{2t}V^{1/2})$.

Now, we use the martingale approximation of $\mathbf{N}(n)-n \bolds
{\rho}$
and $\hat{\bolds{\theta}}(n)-\bolds{\theta}$ from Hu and Zhang
(\citeyear{HZ04}). Let
$\mathcal{F}_m =\sigma(\mathbf{T}_1, \ldots, \mathbf{T}_m,
\mathbf{X}_1,\ldots,\mathbf{X}_m)$ be the $\sigma$-field generated
by the
previous $m$ stages. Then under $\mathcal{F}_{m-1}$, $\mathbf{T}
_m$ and
$\mathbf{X}_m$ are independent, and
\[
E[T_{m1}|F_{m-1}]=g \biggl(\frac{N_1(m-1)}{m-1},\hat{\rho
}_1(m-1) \biggr).
\]
Let $\mathbf{Q}_n=\sum_{m=1}^n \Delta\mathbf{Q}_m$, where $\Delta
\mathbf{Q}_m=(\Delta\mathbf{Q}_{m,1},\Delta\mathbf
{Q}_{m,2})=(\Delta
Q_{m,1k},\Delta Q_{m,2k}$; $k=1,\ldots, d)$ and $\Delta
Q_{m,jk}=T_{m,j}(X_{m,jk}-\theta_{jk})/\rho_j, j=1,2$. Then
$\mathbf{Q}_n=\break O(\sqrt{n\log\log n})$ a.s. is a sequence of martingales
and we can prove
%
%e5.3 ###
\begin{equation}
\hat{\bolds{\theta}}(n)-\bolds{\theta}=\frac{\mathbf
{Q}_n}{n}+O \biggl(\frac
{\log\log
n}{n} \biggr)\qquad\mbox{a.s.}
\end{equation}
Let $\mathbf{M}_n=\sum_{m=1}^n \Delta\mathbf{M}_m$, where $\Delta
\mathbf{M}_m=\mathbf{T}_m-E[\mathbf{T}_m|\mathcal{F}_{m-1}]$, and
$B_{n,m}$ as
defined in Hu and Zhang (\citeyear{HZ04}), then
\begin{eqnarray*}
\mathbf{N}(n)-n\bolds{\rho}&=&\sum_{m=1}^n \Delta\mathbf{M}_m
B_{n,m}+\sum
_{m=1}^n \Delta\mathbf{Q}_m
\nabla(\bolds{\rho})\Big|_{\bolds{\theta}} E \sum_{k=m}^n \frac{1}{k}
B_{n,k}+o(n^{-1/2-\delta/3})
\\
:\!&=&\mathbf{U}_n +o (n^{-1/2-\delta/3} )
\end{eqnarray*}
almost surely, where $\mathbf{U}_n$ is a sum of martingale differences.

We can approximate the process $\mathbf{N}([nt])-[nt]\bolds{\rho}$ and
$\hat{\bolds{\theta}}([nt])-\bolds{\theta}$ (for any point $t \in(0,1]$)
similarly as $\mathbf{N}(n)-n\bolds{\rho}$ and
$\hat{\bolds{\theta}}(n)-\bolds{\theta}$. We obtain
%
%e5.4 ###
\begin{equation}
\hat{\bolds{\theta}}([nt])-\bolds{\theta}=\frac{\mathbf{Q}
_{[nt]}}{[nt]}+O \biggl(\frac{\log\log
[nt]}{[nt]} \biggr)\qquad\mbox{a.s.}
\end{equation}
and
\begin{eqnarray*}
&&\mathbf{N}([nt])-[nt]\bolds{\rho}\\
&&\qquad\!\phantom{:}=\sum_{m=1}^{[nt]} \Delta\mathbf{M}_m
B_{[nt],m}+\sum_{m=1}^{[nt]} \Delta\mathbf{Q}_m
\nabla(\bolds{\rho})\Big|_{\bolds{\theta}} E \sum_{k=m}^{[nt]} \frac{1}{k}
B_{[nt],k}
+o ( ([nt] )^{-1/2-\delta/3} )\\
&&\qquad:=\mathbf{U}_{[nt]} +o (([nt])^{-1/2-\delta/3} )
\end{eqnarray*}
almost
surely.

Hu and Zhang (\citeyear{HZ04}) proved that both martingales $\mathbf{Q}_n$ and
$\mathbf{U}_n$ satisfy the Lindberg conditions. Similarly, we can show
that both
martingales $\mathbf{Q}_{[nt]}$ and $\mathbf{U}_{[nt]}$ also satisfy the
Lindberg conditions. Now we just have to calculate the covariance
matrix of the martingales $\mathbf{Q}_{[nt]}$ and $\mathbf{U}_{[nt]}$.
First, based on the results of Hu and Zhang (\citeyear{HZ04}), we have
\[
\hat{\bolds{\rho}}(n)-\bolds{\rho}=O \Biggl(\sqrt{\frac{\log\log
n}{n}} \Biggr)
\quad\mbox{and}\quad n^{-1}\mathbf{N}(n)-\bolds{\rho}=O \Biggl(\sqrt{\frac
{\log
\log
n}{n}} \Biggr)
\]
almost surely. Therefore, for any $t \in(0,1]$, we
have
%
%e5.5 ###
\begin{eqnarray} \label{consist1}
\hat{\bolds{\rho}}([nt])-\bolds{\rho}&=&O \Biggl(\sqrt{\frac{\log
\log
[nt]}{[nt]}} \Biggr) \quad\mbox{and }\nonumber\\[-8pt]\\[-8pt]
([nt])^{-1}\mathbf{N}([nt])-\bolds{\rho}&=&O \Biggl(\sqrt{\frac{\log
\log
[nt]}{[nt]}} \Biggr)\nonumber
\end{eqnarray}
almost surely. Now, we can calculate $\operatorname{Var} [\Delta
\mathbf{M}_{[nt]}|\mathcal{F}_{[nt]-1} ]$, $\operatorname{Var} [\Delta
\mathbf{Q}_{[nt]}|\mathcal{F}_{[nt]-1} ]$ and $\operatorname{Cov} [\Delta
\mathbf{M}_{[nt]}, \Delta\mathbf{Q}_{[nt]}|\mathcal
{F}_{[nt]-1} ]$.

First, $\Delta\mathbf{M}_{[nt]}=\mathbf{T}_{[nt]}-E [\mathbf{T}
_{[nt]}|\mathcal{F}_{[nt]-1} ]$ is a binary random vector. Based
on conditions
(A2), (A3) and (\ref{consist1}), we have
%
%e5.6 ###
\begin{equation} \label{consist2}
\operatorname{Var} \bigl[\Delta\mathbf{M}_{[nt]}|\mathcal{F}_{[nt]-1}
\bigr]=\Sigma_1 +
o(1)
\end{equation}
almost surely. Similarly, we can show
%
%e5.7 ###
\begin{equation} \label{consist3}
\operatorname{Var} \bigl[\Delta\mathbf{Q}_{[nt]}|\mathcal{F}_{[nt]-1} \bigr]=V+o(1)
\end{equation}
and
%
%e5.8 ###
\begin{equation} \label{consist4}
\operatorname{Cov} \bigl[\Delta\mathbf{M}_{[nt]}, \Delta\mathbf{Q}_{[nt]}|\mathcal
{F}_{[nt]-1} \bigr]=o(1)
\end{equation}
almost surely.

Based on results (\ref{consist2}), (\ref{consist3}) and
(\ref{consist4}), it follows that for any $0<s<t<1$,
\begin{eqnarray*}
\operatorname{Cov} \bigl[\mathbf{Q}_{[ns]},\mathbf{Q}_{[nt]} \bigr]&=&\operatorname{Cov}\Biggl(\sum
_{m=1}^{[ns]}\Delta\mathbf{Q}_m,
\sum_{m=1}^{[nt]}\Delta\mathbf{Q}_m \Biggr)\\
&=&ns\bigl(V+o(1)\bigr)=nsV+o(n),
\\
\operatorname{Cov} \bigl[\mathbf{U}_{[ns]},\mathbf{U}_{[nt]} \bigr]&=&n\wedge
_{11}(s,t)+o(n),
\\
\operatorname{Cov} \bigl[\mathbf{Q}_{[ns]},\mathbf{U}_{[nt]} \bigr]&=&\operatorname{Cov} \Biggl[\sum
_{m=1}^{[ns]}\Delta
\mathbf{Q}_m, \sum_{m=1}^{[nt]}\Delta\mathbf{M}_m
B_{[nt],m}\\
&&\qquad\hspace*{48.9pt}{} + \sum_{m=1}^{[nt]}\Delta\mathbf{Q}_m
\nabla(\bolds{\rho})\Big|_{\bolds{\theta}} E\sum_{k=m}^{[nt]}\frac
{1}{k}B_{[nt],k} \Biggr]\\
&=&\operatorname{Cov} \Biggl[\sum_{m=1}^{[ns]}\Delta\mathbf{Q}_m,\sum
_{m=1}^{[nt]}\Delta
\mathbf{M}_m
B_{[nt],m} \Biggr]\\
&&{} +\operatorname{Cov} \Biggl[\sum_{m=1}^{[ns]}\Delta\mathbf{Q}_m,\sum
_{m=1}^{[nt]}\Delta
\mathbf{Q}_m \nabla(\bolds{\rho})\Big|_{\bolds{\theta}} E\sum
_{k=m}^{[nt]}\frac
{1}{k}B_{[nt],k} \Biggr]\\
&=&\operatorname{Cov} \Biggl[\sum_{m=1}^{[ns]}\Delta\mathbf{Q}_m,\sum
_{m=1}^{[nt]}\Delta
\mathbf{Q}_m \nabla
(\bolds{\rho})\Big|_{\bolds{\theta}} E\sum_{k=m}^{[nt]}\frac
{1}{k}B_{[nt],k} \Biggr]\\
&=&\bigl(V\nabla
(\bolds{\rho})|_{\bolds{\theta}} E+o(1)\bigr)\sum_{m=1}^{[ns]}
\Biggl(\sum
_{k=m}^{[nt]}\frac{1}{k}B_{[nt],k} \Biggr)\\
&=&nV\nabla(\bolds{\rho})\big|_{\bolds{\theta}} E \int_0^s dx
\biggl[\int_x^t
\frac{t}{y} \biggl(\frac{t}{y} \biggr)^H \,dy \biggr]+o(n)\\
&=&n\wedge_{21}(s,t)+o(n),
\end{eqnarray*}
and similarly,
\[
\operatorname{Cov} \bigl[\mathbf{Q}_{[nt]},\mathbf{U}_{[ns]} \bigr]
=n\wedge_{12}(s)+o(n),
\]
where
\begin{eqnarray*}
\wedge_{11}(s,t)&=&\int_0^s  \biggl(\frac{s}{x} \biggr)^{H'} \Sigma_1
\biggl(\frac{t}{x} \biggr)^{H}\,dx\\
&&{} + \int_0^s dx \biggl[\int_x^s
\frac{1}{y} \biggl(\frac{s}{y} \biggr)^H \,dy \biggr]' \Sigma_2
\biggl[\int_x^t \frac{1}{y} \biggl(\frac{t}{y} \biggr)^H \,dy \biggr],
\\
\wedge_{21}(s,t)&=&V\nabla(\bolds{\rho})\big|_{\bolds{\theta}} E \int
_0^s dx
\biggl[\int_x^t\frac{t}{y} \biggl(\frac{t}{y} \biggr)^H \,dy \biggr],
\\
\wedge_{12}(s)&=&\int_0^s dx  \biggl[\int_x^s \frac{s}{y}
\biggl(\frac{s}{y} \biggr)^H \,dy \biggr]E'\nabla(\bolds{\rho}
)\big|_{\bolds{\theta}} 'V.
\end{eqnarray*}

Therefore, the asymptotic covariance function of
$n^{-1/2}(\mathbf{U}_{[nt]},\mathbf{Q}_{[nt]})$ agrees with that of
$(G_t,W_{2t} V^{1/2})$. So by weak convergence of the martingale
[cf. Theorem 4.1 of Hall and Heyde (\citeyear{HH80})], we have
\[
n^{-1/2}([nt]) \biggl(\frac{\mathbf{N}([nt])}{[nt]}-\bolds{\rho},
\hat
{\bolds{\theta}}([nt])-\bolds{\theta} \biggr)
\rightarrow(G_t,W_{2t}V^{1/2})
\]
in distribution in the space
$D_{[0,1]}$ with the Skorohod topology.
\end{pf}
\begin{pf*}{Proof of Theorem \ref{theorem21}}
We assume for $j=1,2$
\[
[nt]\widehat{\operatorname{Var}} (h (\hat{\bolds{\theta}}_j([nt])
) )=
[nt]v_j \bigl(\mathbf{N}([nt])/[nt],\hat{\bolds{\theta
}}([nt]) \bigr)\bigl(1+o_P(1)\bigr)
\]
and
\[
[nt]\operatorname{Var} (h(\hat{\bolds{\theta}}_j([nt])) )=[nt]v_j(\bolds
{\rho}
,\bolds{\theta}),
\]
where $v$ is a continuous function. We also assume
\[
[nt]v_j \bigl(\mathbf{N}([nt])/[nt],\hat{\bolds{\theta
}}([nt]) \bigr)
=[nt]v_j(\bolds{\rho},\bolds{\theta})+O \Biggl(\sqrt{\frac{\log
\log
[nt]}{[nt]}} \Biggr) \qquad\mbox{a.s.},
\]
which holds for most circumstances, since
\[
\mathbf{N}([nt])/[nt]=\bolds{\rho}+O \Biggl(\sqrt{\frac{\log\log
[nt]}{[nt]}} \Biggr) \qquad\mbox{a.s.}
\]
and
\[
\hat{\bolds{\theta}}([nt])=\bolds{\theta}+O \Biggl(\sqrt{\frac
{\log\log
[nt]}{[nt]}} \Biggr) \qquad\mbox{a.s.}
\]
So
\[
[nt]\widehat{\operatorname{Var}} (h (\hat{\bolds{\theta}}_j([nt])
) )
=[nt]\operatorname{Var} (h (\hat{\bolds{\theta}}_j([nt]) )
)+O_P \Biggl(\sqrt{\frac{\log\log[nt]}{[nt]}} \Biggr).
\]
That is,
$[nt]\widehat{\operatorname{Var}} (h (\hat{\bolds{\theta}}_j([nt])
) )$
converges to
$[nt]\operatorname{Var} (h (\hat{\bolds{\theta}}_j([nt]) ) ),j=1,2,$
in probability. By Slutsky's theorem, the sequential statistics
\[
B_t \biggl(\frac{\mathbf{N}([nt])}{[nt]},\hat{\bolds{\theta
}}([nt])
\biggr)=\sqrt{t}
\frac{h (\hat{\bolds{\theta}}_1([nt]) )-h (\hat
{\bolds{\theta}}_2([nt]) )}
{\sqrt{\widehat{\operatorname{Var}} (h (\hat{\bolds{\theta}}_1([nt])
) )
+\widehat{\operatorname{Var}} (h (\hat{\bolds{\theta}}_2([nt]) ) )}}
\]
and
\[
B_t^* (\hat{\bolds{\theta}}([nt]) )=\sqrt{t}
\frac{h (\hat{\bolds{\theta}}_1([nt]) )-h (\hat
{\bolds{\theta}}_2([nt]) )}
{\sqrt{\operatorname{Var} (h (\hat{\bolds{\theta}}_1([nt]) )
)+\operatorname{Var} (h (\hat{\bolds{\theta}}_2([nt]) ) )}}
\]
have the same distribution asymptotically. So we only need to prove
the sequential statistics $B_t^*$ converges to Brownian motion in
distribution. Now
\begin{eqnarray*}
h (\hat{\bolds{\theta}}_j )-h(\bolds{\theta}_j)&=&
(\hat
{\bolds{\theta}}_j-\bolds{\theta}_j )
\bigl(\partial h(\bolds{\theta}_j)/\partial\bolds{\theta}_j\bigr)'+o (\|
\hat
{\bolds{\theta}}_j-\bolds{\theta}_j\|^{1+\delta} )\\
&=& (\hat{\bolds{\theta}}_j-\bolds{\theta}_j )\bigl(\partial
h(\bolds{\theta}_j)/\partial
\bolds{\theta}_j\bigr)'+o (n^{-1/2-\delta/3} ) \qquad\mbox{a.s.},
j=1,2.
\end{eqnarray*}
It is easy to see that
\[
\operatorname{Var} [\hat{\bolds{\theta}}_j([nt]) ]=\operatorname{Var} [\hat
{\bolds{\theta}
}_j(n) ]/t + o(n^{-1}) \qquad\mbox{a.s.},j=1,2.
\]
Here, we define
\[
C=\sqrt{\operatorname{Var} [h(\hat{\bolds{\theta}}_1([nt])) ]+
\operatorname{Var} [h(\hat{\bolds{\theta}}_2([nt]))
]}\sqrt{\operatorname{Var}
[h(\hat
{\bolds{\theta}}_1([ns])) ]
+\operatorname{Var} [h(\hat{\bolds{\theta}}_2([ns])) ]}
\]
and
\[
\mathbf{D}=\bigl(\partial
h(\bolds{\theta}_1)/\partial\bolds{\theta}_1,-\partial
h(\bolds{\theta}_2)/\partial\bolds{\theta}_2\bigr).
\]
Then
\begin{eqnarray*}
B_t^* (\hat{\bolds{\theta}}([nt]) )&=&\sqrt{t}
\frac{h (\hat{\bolds{\theta}}_1([nt]) )-h (\hat
{\bolds{\theta}}_2([nt]) )}
{\sqrt{\operatorname{Var} (h (\hat{\bolds{\theta}}_1([nt]) )
)+\operatorname{Var} (h (\hat{\bolds{\theta}}_2([nt]) ) )}}\\
&=&\sqrt{t}
\frac{h(\bolds{\theta}_1)-h(\bolds{\theta}_2)+ (\hat{\bolds
{\theta}
}([nt])-\bolds{\theta} )
\mathbf{D}'+o (n^{-1/2-\delta/3} )}
{\sqrt{\operatorname{Var} (h (\hat{\bolds{\theta}}_1([nt]) )
)+\operatorname{Var} (h (\hat{\bolds{\theta}}_2([nt]) ) )}}.
\end{eqnarray*}
By the conclusion of Lemma \ref{lemma1}:
\[
n^{-1/2}([nt]) \bigl(\hat{\bolds{\theta}}([nt])-\bolds{\theta
} \bigr)
\rightarrow(W_{2t}V^{1/2})
\]
in distribution in the space $D_{[0,1]}$
with the Skorohod topology. Under $H_0$, we have
\[
B_t^*= \sqrt{t}\frac{ n^{1/2}([nt])^{-1}W_{2t}V^{1/2}\mathbf{D}'}
{\sqrt{\operatorname{Var} (h (\hat{\bolds{\theta}}_1([nt]) ) )
+\operatorname{Var} (h (\hat{\bolds{\theta}}_2([nt]) )
)}}+o
(n^{-\delta/3} )
\]
almost surely. So the sequential statistics $B_t^*$ converges to a
Gaussian process in distribution. In order to prove that $B_t^*$
converges to a ``Brownian process'' in distribution, it is enough to
show $E B_t^* \rightarrow0$ and for any $0<s<t<1$,
\begin{eqnarray*}
\operatorname{cov}(B_t^*,B_s^*) &\rightarrow& s,
\\
\operatorname{cov} (B_t^*,B_s^* )&=&\frac{n\sqrt{ts}}{[nt][ns]} \frac
{\operatorname{cov}(W_{2t},W_{2s})\mathbf{D}V \mathbf{D}'}{C}
+o (n^{-\delta/3} )\\
&=&\frac{n\sqrt{t}s^{3/2}}{[nt][ns]} \frac{\mathbf{D}V \mathbf{D}
'}{C}+o (n^{-\delta/3} )\\
&=&\frac{n\sqrt{t}s^{3/2}}{[nt][ns]}\frac{\mathbf{D} (n
\operatorname{Var}
[\hat{\bolds{\theta}}(n)
-\bolds{\theta} ]+o(1) )\mathbf{D}'}{C}+o
(n^{-\delta
/3} )\\
&=&\frac{n^2\sqrt{t}s^{3/2}}{[nt][ns]} \frac{\partial
h(\bolds{\theta}_1)/\partial
\bolds{\theta}_1\operatorname{Var} [\hat{\bolds{\theta}}_1(n)-\bolds{\theta
}_1
]\,\partial
h(\bolds{\theta}_1)/\partial\bolds{\theta}_1'
}{C}\\
&&{}+\frac{n^2\sqrt{t}s^{3/2}}{[nt][ns]} \frac{\partial
h(\bolds{\theta}_2)/\partial
\bolds{\theta}_2\operatorname{Var} [\hat{\bolds{\theta}}_2(n)-\bolds{\theta
}_2
]\,\partial
h(\bolds{\theta}_2)/\partial\bolds{\theta}_2'}{C}+o(1)\\
&=&\frac{n^2\sqrt{t}s^{3/2}}{[nt][ns]}
\frac{\operatorname{Var} [h(\hat{\bolds{\theta}}_1(n))
]+\operatorname{Var}
[h(\hat
{\bolds{\theta}}_2(n)) ]}
{C}+o(1)\\
&=&\frac{n^2ts^2}{[nt][ns]}+o(1)\\
& \rightarrow & s \qquad\mbox{a.s.}
\end{eqnarray*}
It is easy to see that $E B_t^* \rightarrow0$. This completes the
proof and shows that
$B_t$ is asymptotical Brownian process in distribution.

Under $H_1,$ the sequential statistics
\begin{eqnarray*}
B_t^* (\hat{\bolds{\theta}}([nt]) )&=&\sqrt{t}
\frac{h (\hat{\bolds{\theta}}_1([nt]) )-h (\hat
{\bolds{\theta}}_2([nt]) )
-(h (\bolds{\theta}_1 )-h (\bolds{\theta}_2 ))}
{\sqrt{\widehat{\operatorname{Var}} (h (\hat{\bolds{\theta}}_1([nt])
) )
+\widehat{\operatorname{Var}} (h (\hat{\bolds{\theta}}_2([nt]) )
)}}\\
&&{}+\sqrt{t}
\frac{h (\bolds{\theta}_1 )-h (\bolds{\theta
}_2 )}
{\sqrt{\widehat{\operatorname{Var}} (h (\hat{\bolds{\theta}}_1([nt])
) )
+\widehat{\operatorname{Var}} (h (\hat{\bolds{\theta}}_2([nt]) ) )}}.
\end{eqnarray*}
With similar proof, the first term converges to a standard Brownian
motion in distribution asymptotically.
Because
\[
[nt]\widehat{\operatorname{Var}} (h (\hat{\bolds{\theta}}_j ([nt]
) ) )=
v_j \biggl(\frac{\mathbf{N} ([nt] )}{[nt]},\hat{\bolds
{\theta}
} ([nt] ) \biggr)\bigl(1+o(1)\bigr)
\qquad\mbox{a.s. } j=1,2,
\]
we have that
\[
\sqrt{t}
\frac{h (\bolds{\theta}_1 )-h (\bolds{\theta
}_2 )}
{\sqrt{\widehat{\operatorname{Var}} (h (\hat{\bolds{\theta}}_1([nt])
) )
+\widehat{\operatorname{Var}} (h (\hat{\bolds{\theta}}_2([nt]) ) )}}
\]
converges to
%
%e5.9 ###
\begin{equation}
t\frac{\sqrt{n}(h (\bolds{\theta}_1 )-h (\bolds
{\theta}
_2 ))}
{\sqrt{v_1(\bolds{\rho},\bolds{\theta})+v_2(\bolds{\rho},\bolds
{\theta})}}=\sqrt
{n}\mu t
\end{equation}
in probability. Therefore, under $H_1$, by Slutsky's theorem,
$B_t^*-\sqrt{n}\mu t$ converges to
a standard Brownian motion asymptotically.
\end{pf*}
\end{appendix}

\section*{Acknowledgments}
Special thanks go to anonymous referees, the Associate Editor and the
Editor for the constructive comments, which led to a much improved
version of the paper.

\printaddresses

\end{document}